\documentclass[10pt,french]{amsart}
\usepackage{helvet}
\usepackage[T1]{fontenc}
\usepackage[latin1]{inputenc}
\usepackage{a4wide}
\usepackage{amssymb,latexsym,amsfonts}

 \theoremstyle{plain}    
 \newtheorem{thm}{Théorème}[section]
 \theoremstyle{definition}
  \newtheorem{example}[thm]{Exemple}
 \theoremstyle{definition}
 \newtheorem{defn}[thm]{Définition}
 \theoremstyle{remark}
 \newtheorem{rem}[thm]{Remarque}
 \theoremstyle{plain}    
 \newtheorem{prop}[thm]{Proposition} 
 \theoremstyle{plain}    
 \newtheorem{cor}[thm]{Corollaire} 

\theoremstyle{definition}
\newtheorem{enavant}[thm]{}

\usepackage{babel}
\begin{document}

\title{Quelques remarques sur la notion de modification affine}

 \insert\footins{\footnotesize\noindent \textbf{Mathematics Subject Classification (2000)}: 14E05,14Rxx.

\noindent \textbf{Key words}:  Affine modifications. \unskip\strut\par}

\author{{\bf A. DUBOULOZ}} 
\address{ A. DUBOULOZ 
 \newline\indent INSTITUT FOURIER,
 \newline\indent UMR5582 (UJF-CNRS) 
 \newline\indent BP 74, 38402 ST MARTIN D'HERES CEDEX \\
 \newline\indent  FRANCE}
\email {adrien.dubouloz@ujf-grenoble.fr} 

\begin{abstract}
La notion de modification affine introduite par Kaliman et Zaidenberg
\cite{KaZa99} permet d'investiguer de manière précise la structure
des morphismes birationnels entre variétés affines. Dans cette note,
on construit un contrepartie globale de cette notion. On en déduit
une description analogue de la structure des morphismes affines birationnels
entre variétés quasi-projectives.

\indent\newline\noindent \textsc{\tiny ABSTRACT}. We construct a
global counterpart to the notion of affine modification due to Kaliman
and Zaidenberg \cite{KaZa99}. This leads to an explicit description
of the structure of birational affine morphisms between arbitrary
quasi-projective varieties.
\end{abstract}
\maketitle
\section*{Introduction}

Il est bien connu que tout morphisme projectif birationnel $\pi:Y\rightarrow X$
entre variétés quasi-projectives sur un corps $k$ s'identifie à l'éclatement
d'un sous-schéma fermé de $X$. Il est naturel de chercher une description
analogue de la structure des morphismes birationnels $\pi:Y\rightarrow X$
entre variétés affines, et plus généralement des morphismes affines
birationnels entre variétés quelconques. On peut remarquer tout de
suite que si $\pi:Y\rightarrow X$ est un morphisme affine birationnel,
alors il existe une immersion ouverte $i:Y\hookrightarrow\bar{Y}$
de $Y$ dans un $X$-schéma projectif $\bar{\pi}:\bar{Y}\rightarrow X$
birationnel à $X$. On peut ainsi identifier $Y$ à un ouvert d'un
schéma obtenu à partir de $X$ en faisant éclater un certain sous-schéma
fermé $Z$ de $X$. Si l'on suit cette approche, il reste ensuite à caractériser précisément la structure du bord $\bar{Y}\setminus Y$.

Lorsque $X$ est une variété affine, Kaliman et Zaidenberg \cite{KaZa99}
construisent explicitement un $X$-schéma projectif $\bar{\pi}:\bar{Y}\rightarrow X$
pour lequel $Y$ s'identifie au complémentaire de la \emph{transformé
propre} (cf. (\ref{txt-def:Transforme_propre}) ci-dessous) d'un diviseur
principal $D=\textrm{div}\left(f\right)$ sur $X$ dont le support
contient le sous-schéma fermé $Z$ que l'on fait éclater. 
Cela conduit à la notion de \emph{modification affine}, qui englobe 
certaines constructions antérieures introduites entre autres par Zariski
\cite{Zar43} et Davis \cite{Dav67}. 

Notre but ici est de généraliser le procédé de Kaliman et Zaidenberg
afin de décrire la structure des morphismes affines birationnels $\pi:Y\rightarrow X$
entre variétés quasi-projectives. Idéalement, on cherche à construire
un sous-schéma fermé $Z$ de $X$, contenu dans un diviseur de Cartier
effectif $D$, tel que $Y$ s'identifie au complémentaire dans l'éclatement
$\bar{\pi}:\bar{Y}\rightarrow X$ de $Z$ de la transformée propre
$D^{pr}$ de $D$, en un sens à préciser. Une approche naïve du problème
consisterait simplement à appliquer le procédé précédent sur un recouvrement
affine convenable de $X$ en espérant que les données locales correspondantes
se recollent agréablement. Cela n'est a priori pas le cas, une des
obstructions provenant par exemple du fait que le diviseur principal
$D=\textrm{div}\left(f\right)$ qui apparaît dans la construction
de Kaliman et Zaidenberg n'est pas défini de manière canonique. Notre
méthode consiste par conséquent à définir directement une contrepartie
''globale'' du procédé de Kaliman et Zaidenberg. Cela conduit à
la caractérisation suivante :

\begin{thm}
Soit $X$ une variété projective sur un corps $k$. Pour tout variété
$Y$ et tout morphisme affine birationnel $\pi:Y\rightarrow X$, il
existe un sous-schéma fermé $Z\subset X$, d'idéal $\mathcal{I}\subset\mathcal{O}_{X}$
et un diviseur de Cartier effectif $D\in\textrm{Div}^{+}\left(X\right)$
contenant $Z$, tel que $Y$ soit isomorphe au $X$-schéma \[
\sigma_{\mathcal{I},D}:\tilde{X}_{\mathrm{\mathcal{I},D}}=\mathbf{Spec}\left(\left(\bigoplus_{n\geq0}\left(\mathcal{I}\otimes\mathcal{O}_{X}\left(D\right)\right)^{n}t^{n}\right)/\left(1-t\right)\right)\longrightarrow X,\]
 et tel que $\pi:Y\rightarrow X$ corresponde avec $\sigma_{\mathcal{I},D}$
via cet l'isomorphisme. Le morphisme $\sigma_{\mathcal{I},D}$ ainsi
défini sera appelé la modification affine de $X$ de centre $\left(\mathcal{I},D\right)$. 
\end{thm}
\begin{enavant} Expliquons en quoi l'énoncé ci-dessus constitue une
réponse au problème initial. Tout d'abord, puisque par hypothèse $Z$
est contenu dans $D$, son idéal $\mathcal{I}$ contient l'idéal $\mathcal{O}_{X}\left(-D\right)$
de $D$. Par conséquent, $\mathcal{I}\otimes\mathcal{O}_{X}\left(D\right)$
contient $\mathcal{O}_{X}$, donc en particulier la section constante
$1$, ce qui justifie la définition de $\tilde{X}_{\mathcal{I},D}$.
D'autre part, via l'immersion ouverte canonique \[
j:\tilde{X}_{\mathcal{I},D}\hookrightarrow\mathbf{Proj}_{X}\left(\bigoplus_{n\geq0}\left(\mathcal{I}\otimes\mathcal{O}_{X}\left(D\right)\right)^{n}t^{n}\right)\simeq\mathbf{Proj}_{X}\left(\bigoplus_{n\geq0}\mathcal{I}^{n}t^{n}\right),\]
 $\tilde{X}_{\mathcal{I},D}$ s'identifie au complémentaire dans l'éclatement
$\bar{\pi}:\bar{Y}\rightarrow X$ de $Z$ du sous-schéma fermé \[
D^{pr}=\mathbf{Proj}_{X}\left(\left(\bigoplus_{n\geq0}\left(\mathcal{I}\otimes\mathcal{O}_{X}\left(D\right)\right)^{n}t^{n}\right)/t\right)\]
 de $\bar{Y}$. Cela montre que $\tilde{X}_{\mathcal{I},D}$ est bien
birationnel à $X$, tout en définissant au passage ce qu'est la transformée
propre du diviseur $D$. 

\end{enavant}

\begin{example}
Soit $X$ une variété quasi-projective complexe, $p_{1}:X\times\mathbb{C}\rightarrow X$
le fibré en droites trivial sur $X$, $X_{0}\subset X\times\mathbb{C}$
la section nulle de ce fibré, et soit $\tilde{D}\in\textrm{Div}^{+}\left(X\right)$
un diviseur de Cartier effectif non nul. Soit $Z=p_{1}^{*}\tilde{D}\cap X_{0}$
le sous-schéma fermé de $X\times\mathbb{C}$ défini par l'idéal $\mathcal{I}_{Z}=\mathcal{O}_{X\times\mathbb{C}}\left(-X_{0}\right)+\mathcal{O}_{X\times\mathbb{C}}\left(-p_{1}^{*}\tilde{D}\right)$,
et soit $\sigma_{\mathcal{I}_{Z},p_{1}^{*}\tilde{D}}:L\rightarrow X\times\mathbb{C}$
la modification affine de $X\times\mathbb{C}$ de centre $\left(\mathcal{I}_{Z},D\right)$.
Par construction, \[
L\simeq\mathbf{Spec}_{X\times\mathbb{C}}\left(\bigoplus_{n\geq0}\mathcal{O}_{X\times\mathbb{C}}\left(n\left(p_{1}^{*}\tilde{D}-X_{0}\right)\right)\right)\]
 s'identifie, en tant que $X$-schéma, au fibré en droites $p_{1}\circ\sigma_{\mathcal{I}_{Z},D}:L\rightarrow X$
espace total du faisceau inversible $\mathcal{O}_{X}\left(-\tilde{D}\right)$
sur $X$. Via cet isomorphisme, $\sigma_{\mathcal{I}_{Z},p_{1}^{*}\tilde{D}}$
correspond alors avec le $X$-morphisme induit par la section canonique
de $\mathcal{O}_{X}\left(\tilde{D}\right)$ de diviseur $\tilde{D}$.
\end{example}
Cette note se compose de deux partie. La première passe rapidement
en revue la notion de modification affine selon \cite{KaZa99}. La
seconde est consacrée à la contrepartie ''globale'' de cette notion. 

\nocite{EGAII}

\section{Modifications affines selon Kaliman et Zaidenberg }

Dans ce paragraphe, nous rappelons la notion de modification affine
d'un schéma affine sous une forme légèrement différente de celle introduite
par Kaliman et Zaidenberg \cite{KaZa99}. 

\begin{defn}
\label{def:Modification_affine_sens_large} On appelle \emph{modification
affine} d'un schéma affine $X=\textrm{Spec}\left(A\right)$, tout
$X$-schéma affine de type fini $\pi:Y\rightarrow X$ spectre d'une
sous-$A$-algèbre de type fini de l'anneau total des fractions $K_{A}$
de $A$. 
\end{defn}
\begin{rem}
Lorsque $X$ est intègre, $K_{A}$ s'identifie au corps des fractions
$\textrm{Frac}\left(A\right)$ de $A$. Dans ce cas, toute sous-$A$-algèbre
de type fini $B\subset\textrm{Frac}\left(A\right)$ définit un schéma
affine $Y=\textrm{Spec}\left(B\right)$. Puisque $A$ et $B$ ont
le même corps de fractions, l'inclusion $A\hookrightarrow B$ induit
un morphisme birationnel $\pi:Y\rightarrow X$. 
\end{rem}
\begin{enavant} \label{txt:Presention_KaZa} Une sous-$A$-algèbre
de type fini de $B\subset K_{A}$ est engendrée par un sous-$A$-module
de type fini $J$ de $K_{A}$, encore appelé \emph{idéal fractionnaire.}
L'ensemble $F=\left\{ f\in A,\; fJ\subset A\right\} $ des \emph{dénominateurs
de} $J$ est un idéal de $A$. Le choix d'un élément $f\in F\cap K_{A}^{*}$
permet alors identifier $J$ au sous-$A$-module $f^{-1}\tilde{J}\subset A_{f}$,
où $\tilde{J}=fJ\subset A$. En posant $I=\left(f,\tilde{J}\right)\subset A$,
on obtient finalement des isomorphismes de $A$-algèbres \[
B\simeq A\left[J\right]\simeq A\left[I/f\right]\subset A_{f}\subset K_{A},\]
 où $A\left[I/f\right]=\left\{ a^{k}/f^{k}\textrm{ }\mid\textrm{ }a^{k}\in I^{k},\textrm{ }k\geq1\right\} $
s'identifie canoniquement au quotient $A\left[\left(f^{-1}I\right)t\right]/\left(1-t\right)$
de la $A$-algèbre de Rees \[
A\left[\left(f^{-1}I\right)t\right]=\bigoplus_{n\geq0}\left(f^{-1}I\right)^{n}t^{n}\]
 du sous-$A$-module $f^{-1}I=\left(1,J\right)$ de $K_{A}$ par l'idéal
principal $\left(1-t\right)$. 

\end{enavant}

\begin{rem}
La $A$-algèbre $A\left[I/f\right]$ s'identifie aussi canoniquement
au quotient $A\left[It\right]/\left(1-ft\right)$ de la $A$-algèbre
de Rees $A\left[It\right]$ par l'idéal principal $\left(1-ft\right)$.
On retrouve ainsi la définition d'une modification affine selon \cite{KaZa99}.
\end{rem}
\begin{enavant} L'idéal fractionnaire $J$, le dénominateur commun
$f$, ainsi que l'idéal $I$ de $A$ permettant d'obtenir l'identification
précédente ne sont évidemment pas uniques. Nous adopterons néanmoins
la terminologie suivante, issue de \cite{KaZa99}. 

\end{enavant}

\begin{defn}
\label{def:Modif_affine_KaZa} Étant donné un idéal $I$ de $A$ et
un élément régulier $f\in I$, le $X$-schéma \[
\sigma_{I/f}:X_{I,f}=\textrm{Spec}\left(A\left[I/f\right]\right)\rightarrow X=\textrm{Spec}\left(A\right)\]
est appelé \emph{modification affine de} $X$ \emph{de centre} $\left(I,f\right)$. 
\end{defn}
\begin{enavant} \label{txt-def:Transforme_propre} Puisque $A\left[I/f\right]\simeq A\left[It\right]/\left(1-ft\right)$,
l'immersion canonique $ 
j:X_{I/f}\hookrightarrow\tilde{X}=\textrm{Proj}_{A}\left(A\left[It\right]\right)$
 réalise $X_{I,f}$ comme l'ouvert principal $\tilde{X}_{f}$ de $\tilde{X}$
où $f$ ne s'annule pas. En conséquence, le morphisme $\sigma:X_{I/f}\rightarrow X$
se factorise, via l'immersion $j$, par l'éclatement $\textrm{Bl}_{I}:\tilde{X}\rightarrow X$
de $X$ de centre $I$. Suivant \cite{KaZa99}, nous dirons que le
sous-schéma fermé $D_{f}^{pr}\simeq\textrm{Proj}_{A}\left(A\left[It\right]/ft\right)$,
complémentaire dans $\tilde{X}$ de $j\left(X_{\mathcal{I},f}\right)\simeq\tilde{X}_{f}$
est \emph{la transformée propre} de $D$. L'exemple suivant montre
que la transformée propre $D_{f}^{pr}$ est en général distincte de
la \emph{transformée stricte} $D_{f}'$ de $D=\textrm{div}\left(f\right)$,
définie comme l'adhérence dans $\tilde{X}$ de l'image inverse de
$D\setminus V\left(I\right)$ par l'éclatement $\textrm{Bl}_{I}:\tilde{X}\rightarrow X$. 

\end{enavant}

\begin{example}
Soient $k$ un corps, $A=k\left[x,y\right]$, $I=\left(x,y\right)$
et $f=x^{2}\in I$. La transformée propre du diviseur $\textrm{div}\left(f\right)$
dans $\textrm{Proj}_{A}A\left[xt,yt\right]\simeq\textrm{Proj}_{A}A\left[u,v\right]/\left(xv-yu\right)$
est donnée par l'équation $\left\{ xu=0\right\} $. Elle est donc
constituée de la réunion de la transformée stricte de la droite $V\left(f\right)=\left\{ x=0\right\} $
et du diviseur exceptionnel $\left\{ x=y=0\right\} $ de l'éclatement
$\textrm{Bl}_{I}:\textrm{Proj}_{A}A\left[xt,yt\right]\rightarrow\textrm{Spec}\left(A\right)$.
\end{example}
\begin{enavant} Il existe cependant des cas où la transformée propre
et la transformée stricte de $D$ coïncident, comme le montre la proposition
suivante (cf. aussi \cite{Dav67}) : 

\end{enavant}

\begin{prop}
\label{pro:Suite_regulier_implique_TPr_egal_TS} Si l'idéal $I\subset A$
est engendré par une suite régulière $a_{0}=f,a_{1},\ldots,a_{r}$,
alors la transformée propre de $D=\textrm{div}\left(f\right)$ coïncide
avec sa transformée stricte.
\end{prop}
\begin{proof}
Puisque la suite $a_{0},\ldots,a_{r}$ est régulière, l'algèbre de
Rees $A\left[It\right]$ s'identifie, via l'homomorphisme $X_{i}\mapsto a_{i}t$,
$i=0,\ldots,r$, à la $A$-algèbre $B$, quotient de $A\left[X_{0},\ldots,X_{r}\right]$
par l'idéal engendré par les mineurs $2\times2$ de la matrice \[
\left(\begin{array}{ccc}
a_{0} & \cdots & a_{r}\\
X_{0} & \cdots & X_{r}\end{array}\right).\]
 Ainsi la transformée propre $D_{f}^{pr}$ et la transformée stricte
$D_{f}'$ coïncident puisqu'elles sont toutes deux données par l'équation
$\left\{ X_{0}=0\right\} $ dans $\textrm{Proj}_{A}B$. 
\end{proof}

\section{Modifications affines générales }

D'après ce qui précède, tout morphisme birationnel $\sigma:Y\rightarrow X$
entre variétés affines s'identifie à une modification affine de $X$
de centre $\left(I,f\right)$ convenable. Dans ce paragraphe, nous
cherchons à construire une contrepartie ''globale'' de ce résultat,
à savoir une description de la structure des morphismes affines birationnels
entre schémas quelconques en termes analogues aux précédents.

\subsection{Morphismes affines birationnels}

Rappelons tout d'abord qu'étant donné un schéma $X$ fixé, on dit
qu'un $X$-schéma $\pi:Y\rightarrow X$ est affine sur $X$ s'il existe
un recouvrement $\left(X_{\alpha}\right)$ de $X$ par des ouverts
affines, tel que pour tout $\alpha$, $\pi^{-1}\left(X_{\alpha}\right)$
soit un ouvert affine de $Y$. 

\begin{enavant} Dans ce contexte ''global'', la correspondance
bi-univoque entre les schémas affines de type fini $X$ sur un corps
$k$ et les $k$-algèbres de type fini $B=\Gamma\left(X,\mathcal{O}_{X}\right)$,
se transpose en une correspondance entre les $X$-schémas affines
de type fini $\pi:Y\rightarrow X$ et les $\mathcal{O}_{X}$-algèbres
quasi-cohérentes de type fini $\mathcal{B}=\pi_{*}\mathcal{O}_{Y}$.
Rappelons que le faisceau des germes de fonctions méromorphes $\mathcal{K}_{X}$
d'un schéma $X$ donné est défini comme le faisceau associé au pré-faisceau
défini par $\mathcal{K}_{X}\left(U\right)=K_{A}$ pour tout ouvert
affine $U=\textrm{Spec}\left(A\right)$ de $X$. La notion de sous-$A$-algèbre
de type fini de l'anneau total des fractions $K_{A}$ d'un anneau
$A$ admet alors la contrepartie suivante.

\end{enavant}

\begin{defn}
Soit $X$ un schéma. On appelle $\mathcal{O}_{X}$-\emph{algèbre fractionnaire}
toute $\mathcal{O}_{X}$-algèbre isomorphe à une sous-$\mathcal{O}_{X}$-algèbre
quasi-cohérente de type fini du faisceau des germes de fonctions méromorphes
$\mathcal{K}_{X}$ de $X$. 
\end{defn}
\begin{enavant} Une $\mathcal{O}_{X}$-algèbre fractionnaire $\mathcal{B}\subset\mathcal{K}_{X}$
définit un $X$-schéma affine de type fini $\sigma:Y=\mathbf{Spec}\left(\mathcal{B}\right)\rightarrow X.$
Par construction, l'homomorphisme $\mathcal{O}_{X}\rightarrow\mathcal{B}$
est injectif et induit un isomorphisme $\mathcal{K}_{X}\stackrel{\sim}{\rightarrow}\sigma_{*}\mathcal{K}_{Y}$.
Réciproquement, si $\sigma:Y\rightarrow X$ un morphisme affine de
type fini tel que l'homomorphisme $\sigma^{\sharp}:\mathcal{O}_{X}\rightarrow\sigma_{*}\mathcal{O}_{Y}$
induise un isomorphisme $\tilde{\sigma}:\mathcal{K}_{X}\stackrel{\sim}{\rightarrow}\sigma_{*}\mathcal{K}_{Y}$
alors $\sigma^{\sharp}$ est injectif puisque les homomorphismes canoniques
$i_{X}:\mathcal{O}_{X}\rightarrow\mathcal{K}_{X}$ et $\sigma_{*}\left(i_{Y}\right):\sigma_{*}\mathcal{O}_{Y}\rightarrow\sigma_{*}\mathcal{K}_{Y}$
le sont. Ainsi, puisque $Y$ est de type fini sur $X$, $\mathcal{B}=\tilde{\sigma}^{-1}\left(\sigma_{*}\mathcal{O}_{Y}\right)\subset\mathcal{K}_{X}$
est une $\mathcal{O}_{X}$-algèbre fractionnaire isomorphe à $\sigma_{*}\mathcal{O}_{Y}$.
On obtient donc la caractérisation suivante : 

\end{enavant}

\begin{prop}
Il y a une correspondance biunivoque entre les $\mathcal{O}_{X}$-algèbres
fractionnaires et les morphismes affines de type fini $\pi:Y\rightarrow X$
induisant un isomorphisme entre les faisceaux des germes de fonctions
méromorphes sur $X$ et $Y$ respectivement.
\end{prop}
\noindent Lorsque $X$ est un schéma intègre, toute $\mathcal{O}_{X}$-algèbre
fractionnaire $\mathcal{B}$ est intègre, et le morphisme structural
$\pi:Y=\mathbf{Spec}\left(\mathcal{B}\right)\rightarrow X$ correspondant
est birationnel. 

\begin{cor}
Si $X$ est un schéma intègre, alors il y a une correspondance biunivoque
entre les $\mathcal{O}_{X}$-algèbres fractionnaires et les morphismes
affines birationnels de type fini $\pi:Y\rightarrow X$.
\end{cor}
\begin{example}
Soient $X$ un schéma intègre, $\mathcal{I}\subset\mathcal{O}_{X}$
un idéal de type fini, $\mathcal{S}=\bigoplus_{n\geq0}{\displaystyle \mathcal{I}^{n}}$,
et soit $\sigma:\tilde{X}=\mathbf{Proj}\left(\mathcal{S}\right)\rightarrow X$
l'éclatement de $X$ de centre $\mathcal{I}$. Pour toute section
globale $f\in\Gamma\left(X,\mathcal{I}\right)$, $\mathcal{B}=\mathcal{S}/\left(1-f\right)\mathcal{S}$
est alors une $\mathcal{O}_{X}$-algèbre fractionnaire. Le $X$-schéma
affine correspondant est canoniquement isomorphe à l'ouvert principal
$\tilde{X}_{f}$ de $\tilde{X}$ où $f$ ne s'annule pas, \emph{i.e.}
l'unique ouvert de $\tilde{X}$ tel que pour tout ouvert affine $U\subset X$,
$\tilde{X}_{f}\cap\sigma^{-1}\left(U\right)=D_{+}\left(f\mid_{U}\right)$
dans $\sigma^{-1}\left(U\right)\simeq\textrm{Proj}\left(\Gamma\left(U,\mathcal{S}\right)\right)$. 
\end{example}

\subsection{Modification affine d'un schéma le long d'un diviseur}

Par définition, une $\mathcal{O}_{X}$-algèbre fractionnaire $\mathcal{B}\subset\mathcal{K}_{X}$
est engendrée par idéal fractionnaire \emph{}$\mathcal{J}\subset\mathcal{K}_{X}$
de $\mathcal{K}_{X}$. Cependant, sans hypothèses supplémentaires
sur $X$, le faisceau d'idéaux $\mathcal{F}$ des dénominateurs de
$\mathcal{J}$, \emph{i.e.} le faisceau défini localement par $\left\{ f\in\mathcal{O}_{X}\textrm{ }\mid\textrm{ }f\cdot\mathcal{J}\subset\mathcal{O}_{X}\right\} $,
n'admet plus nécessairement de section globale.

\begin{enavant} Supposons cependant qu'il existe un $\mathcal{O}_{X}$-module
inversible $\mathcal{L}$ et un homomorphisme régulier non nul $f:\mathcal{L}^{-1}\rightarrow\mathcal{F}$.
Puisque $f$ est en particulier injectif, la composition \[
\mathcal{J}\otimes\mathcal{L}^{-1}\stackrel{\textrm{Id}\otimes f}{\longrightarrow}\mathcal{J}\otimes\mathcal{F}\rightarrow\mathcal{J}\cdot\mathcal{F}\subset\mathcal{O}_{X}\]
 permet alors d'identifier $\mathcal{J}$ à un idéal fractionnaire
de la forme $\tilde{\mathcal{J}}\otimes\mathcal{L}=\tilde{\mathcal{J}}\cdot\mathcal{L}\subset\mathcal{K}_{X}$
pour un certain idéal quasi-cohérent de type fini $\tilde{\mathcal{J}}\subset\mathcal{O}_{X}$.
D'autre part, puisque $\mathcal{F}\subset\mathcal{O}_{X}$, l'homomorphisme
$f:\mathcal{L}^{-1}\rightarrow\mathcal{F}\hookrightarrow\mathcal{O}_{X}$
correspond par dualité à une section régulière globale $s$ de $\mathcal{L}$.
En posant $D=\textrm{div}\left(s\right)$, on a donc une identification
$\mathcal{L}\simeq\mathcal{O}_{X}\left(D\right)$, d'où finalement
l'identification $\mathcal{J}\simeq\tilde{\mathcal{J}}\otimes\mathcal{O}_{X}\left(D\right)$. 

\end{enavant}

\begin{enavant} Puisque $D$ un diviseur de Cartier effectif, $\mathcal{L}^{-1}\simeq\mathcal{O}_{X}\left(-D\right)$
est un idéal inversible de $\mathcal{O}_{X}$. Ainsi $\mathcal{I}=\tilde{\mathcal{J}}+\mathcal{O}_{X}\left(-D\right)$
est un idéal quasi-cohérent de type fini de $\mathcal{O}_{X}$, et
l'on a des isomorphismes de $\mathcal{O}_{X}$-algèbres \[
\mathcal{B}\simeq\mathcal{O}_{X}\left[\mathcal{J}\right]\simeq\mathcal{O}_{X}\left[\mathcal{I}\otimes\mathcal{O}_{X}\left(D\right)\right].\]
 On conclut alors comme dans le cas affine que $\mathcal{B}$ s'identifie
à la $\mathcal{O}_{X}$-algèbre fractionnaire $\mathcal{O}_{X}\left[\mathcal{I}/D\right]$,
quotient de la $\mathcal{O}_{X}$-algèbre de Rees\[
\mathcal{O}_{X}\left[\left(\mathcal{I}\otimes\mathcal{O}_{X}\left(D\right)\right)t\right]=\bigoplus_{n\geq0}\left(\mathcal{I}\otimes\mathcal{O}_{X}\left(D\right)\right)^{n}t^{n}\subset\mathcal{K}_{X}\left[t\right]\]
 de l'idéal fractionnaire $\mathcal{I}\otimes\mathcal{O}_{X}\left(D\right)$
par l'idéal principal engendré par $\left(1-t\right)$. Cela conduit
à la définition suivante. 

\end{enavant}

\begin{defn}
Étant donné un triplet $\left(X,\mathcal{I},D\right)$ constitué d'un
schéma $X$, d'un diviseur effectif $D\in\textrm{Div}^{+}\left(X\right)$
et d'un idéal quasi-cohérent de type fini $\mathcal{I}\subset\mathcal{O}_{X}$
contenant l'idéal $\mathcal{O}_{X}\left(-D\right)$, le $X$-schéma
affine \[
\sigma_{\mathcal{I},D}:X_{\mathcal{I},D}=\mathbf{Spec}\left(\mathcal{O}_{X}\left[\mathcal{I}/D\right]\right)\rightarrow X\]
 défini ci-dessus est appelé la \emph{modification affine de} $X$
\emph{de centre} $\left(\mathcal{I},D\right)$. 
\end{defn}
\begin{enavant} Pour des données locales $\left\{ \left(X_{i},f_{i}\right)_{i\in I}\right\} $
du diviseur de Cartier $D$, $\mathcal{O}_{X}\left(D\right)\mid_{X_{i}}$
est isomorphe au sous-$\mathcal{O}_{X_{i}}$-module de $\mathcal{K}_{X_{i}}$
engendré par $f_{i}^{-1}$. La $\mathcal{O}_{X_{i}}$-algèbre $\mathcal{O}_{X}\left[\mathcal{I}/D\right]\mid_{X_{i}}$
est alors isomorphe au quotient de la $\mathcal{O}_{X_{i}}$-algèbre
de Rees de l'idéal $\mathcal{I}_{i}=\mathcal{I}\mid_{X_{i}}$ de $\mathcal{O}_{X_{i}}$
par l'idéal principal engendré par $\left(1-f_{i}t\right)$. Le schéma
$X_{\mathcal{I},D}\mid_{X_{i}}$ correspondant s'identifie, via l'immersion
ouverte canonique \[
j_{i}:X_{\mathcal{I},D}\mid_{X_{i}}\simeq\mathbf{Spec}\left(\mathcal{O}_{X_{i}}\left[\mathcal{I}_{i}t\right]/\left(1-f_{i}t\right)\right)\hookrightarrow\tilde{X}_{i}=\mathbf{Proj}\left(\mathcal{O}_{X_{i}}\left[\mathcal{I}_{i}t\right]\right),\]
 à l'ouvert principal $\left(\tilde{X}_{i}\right)_{f_{i}}$de $\tilde{X}_{i}$.
Lorsque les $X_{i}$ sont affines, on retrouve évidemment la notion
du paragraphe précédent.

\end{enavant}

\begin{example}
\label{exa:(Compl=E9mentaire-d'un-diviseur).} (Complémentaire d'un
diviseur). Pour tout diviseur effectif $D\in\textrm{Div}^{+}\left(X\right)$
non nul, la modification $\sigma_{1,D}:X_{1,D}\rightarrow X$ de centre
$\left(\mathcal{I}=1\cdot\mathcal{O}_{X},D\right)$ s'identifie à
l'immersion ouverte du sous-schéma $X_{s_{D}}$de $\mathbf{Proj}\left(\mathbf{S}\left(\mathcal{O}_{X}\left(D\right)\right)\right)\simeq X$
où la section canonique de $\mathcal{O}_{X}\left(D\right)$ de diviseur
$D$ ne s'annule pas, c'est à dire au complémentaire du support de
$D$ dans $X$.
\end{example}
\noindent On obtient finalement la contrepartie globale suivante
de la notion de modification affine selon \cite{KaZa99}. 

\begin{thm}
Soit $X$ une variété quasi-projective sur un corps $k$. Alors tout
morphisme affine birationnel de type fini $\pi:Y\rightarrow X$ est
une modification affine de $X$ de centre $\left(\mathcal{I},D\right)$
convenable.
\end{thm}
\begin{proof}
Puisque $X$ est intègre, $\mathcal{K}_{X}$ n'est autre que le faisceau
constant $\mathcal{R}_{X}$ des fonctions rationnelles sur $X$. Soit
dont $\mathcal{J}\subset\mathcal{R}_{X}$ un idéal fractionnaire générateur
de $\pi_{*}\mathcal{O}_{Y}$, et soit $\mathcal{F}$ l'idéal des dénominateurs
communs de $\mathcal{J}$. Puisque $X$ est quasiprojectif, il existe
un $\mathcal{O}_{X}$-module inversible ample $\mathcal{L}$. Ainsi,
pour un certain $n\geq0$, $\mathcal{F}\otimes\mathcal{L}^{n}$ admet
un section globale non nulle. On en déduit un homomorphisme injectif
$f:\mathcal{L}^{-n}\rightarrow\mathcal{F}$ à partir duquel on peut
appliquer la construction décrite précédemment. 
\end{proof}
\bibliographystyle{amsplain}

\begin{thebibliography}{1}

\bibitem{Dav67}
E.D. Davis, \emph{Ideals of the principal class, {R}-sequences and a certain
  monoidal transformation}, Pacific J. Math. \textbf{20} (1967), 197--205.

\bibitem{EGAII}
A. Grothendieck and J. Dieudonné, \emph{{EGA II.} {É}tude globale {é}l{é}mentaire de quelques classes de  morphismes}, Publ. Math. IHES, vol. 8, 1961.



\bibitem{KaZa99}
S. Kaliman and M. Zaidenberg, \emph{Affine modifications and affine
  hypersurfaces with a very transitive automorphism group}, Transformation
  Groups \textbf{4} (1999), 53--95.

\bibitem{Zar43}
O. Zariski, \emph{Foundations of a general theory of birational
  correspondences}, Trans. Amer. Math. Soc. \textbf{53} (1943), 497--542.

\end{thebibliography}

\end{document}